# Methods to Find Integer Points on the Elliptic Curve of Factorizable Polynomial Equation

XIAODONG ZHUANG[1,2,3], NIKOS E. MASTORAKIS[1,2]
[1]Technical University of Sofia, Sofia,
BULGARIA

[2]International Research Institute of Biomedicine, Natural Sciences and Technology, Sofia,
BULGARIA

[3]Qingdao University, Qingdao,
CHINA

*Abstract:* - As a valuable theoretical and application problem, the integer points on two types of elliptic curve $y^2=(x-a)(x-b)(x-c)$ and $y^2=(x-a)(x^2+ax+b)$ are studied. By using elementary number theory method, the solution of the original equations is reduced to the solution of simultaneous Pell equations or generalized Pell equation, which can improve the efficiency in searching integer solutions by computer. Effective methods are proposed to find the integer solutions of the equations, which are convenient for algorithm programming implementation. At the same time, the derivation of the methods leads to two sufficient conditions, which can determine the unsolvability of these two types of equation respectively. The methods proposed are for equations with undetermined coefficients, which are more general than the solutions of the equations with specific known coefficients, and more suitable for computer programming implementation.

## 1 Introduction

Finding the integer points on the elliptic curve $y^2=x^3+Ax+B$ belongs to the problem of Diophantine Equation in number theory, which has significant application in cryptography, [1], [2], [3]. Finding the integer solution of the equation is a key problem. Although it has been proved that there is no general algorithm to determine whether any Diophantine equation has integer solution, [4], [5], for a pacific type of Diophantine equation it is still possible to design algorithms to search its integer solution. Baker's method (and its later improvements) can estimate the bound of possible integer solution for some Diophantine equations, which theoretically enables enumeration to search for the integer solution, [6], [7]. But the value range estimated is often quite large, or even computationally impossible with current computers, [8]. For the elliptic curve, there have also been some results of solution bound estimation, but the range is too large as well, [9], [10], [11], [12], [13], [14]. It has practical importance to improve the search efficiency for possible solutions.

In this paper, two specific types of elliptic curve equation: $y^2=(x-a)(x-b)(x-c)$ and $y^2=(x-a)(x^2+ax+b)$ are studied (i.e. the right side of the equation is factorizable in integer-coefficient polynomials). Effective methods are proposed to find their integer solutions by solving generalized Pell equations. And two sufficient conditions of their insolvability are presented, which can be utilized to improve the efficiency of computer algorithm implementation.

## 2 The Solution of $y^2=(x-a)^3$ and $y^2=(x-a)(x-b)^2$

These two cases are trivial. For the equation $y^2=(x-a)^3$, $x$-$a$ must be a square number: $x-a=u^2$. Then the solution is:

$$\begin{cases} x=a+u^2 \\ y=u^3 \end{cases} \quad u\in Z \tag{1}$$

For $y^2=(x-a)(x-b)^2$, we have:

$$\left(\frac{y}{x-b}\right)^2=x-a \tag{2}$$

The left side of (2) is a square number. Therefore, the right-side *x*-*a* must be a square number: $x-a=u^2$. The solution of (2) is:

$$\begin{cases} x = a + u^2 \\ y = u(a + u^2 - b) \end{cases} \quad u \in Z \tag{3}$$

## 3 The Solution of $y^2 = (x-a)(x-b)(x-c)$

Assume $a<b<c$. Let $z=x-a$, and the equation becomes:

$$y^2 = z(z-m)(z-n) \tag{4}$$

where $m=b-a>0$, $n=c-a>0$.

Suppose $d = \gcd(z-m, z-n)$. Then we have $d|((z-m)-(z-n))$, that is, $d|(n-m)$, and also $d|(c-b)$. Therefore, all the possible values of $d$ constitute a subset of the factor set of $(c-b)$, which can be enumerated finitely. Then (4) becomes:

$$y^2 = z \cdot d^2 \cdot \frac{z-m}{d} \cdot \frac{z-n}{d} \tag{5}$$

where $\gcd\left(\frac{z-m}{d} \cdot \frac{z-n}{d}\right) = 1$. Convert (5) to:

$$\left(\frac{y}{d}\right)^2 = z \cdot \frac{z-m}{d} \cdot \frac{z-n}{d} \tag{6}$$

Let $e = \gcd\left(z, \frac{z-m}{d}\right)$, $f = \gcd\left(z, \frac{z-n}{d}\right)$. Then we have $e \mid \left(z - d \cdot \frac{z-m}{d}\right)$, that is, $e|m$. And we also have $f|n$. Therefore, *e* and *f* belong to the factor sets of *m* and *n* respectively, which can be enumerated finitely.

Because $e \mid \left(\frac{z-m}{d}\right)$, $f \mid \left(\frac{z-n}{d}\right)$, $\gcd\left(\frac{z-m}{d}, \frac{z-n}{d}\right) = 1$, it follows that $\gcd(e,f) = 1$. Moreover, because $e|z$, $f|z$, we have $(ef)|\ z$. Therefore, we have:

$$\left(\frac{y}{d}\right)^2 = (ef)^2 \cdot \frac{z}{ef} \cdot \frac{z-m}{de} \cdot \frac{z-n}{df} \tag{7}$$

and:

$$\left(\frac{y}{def}\right)^2 = \left(\frac{z}{ef}\right)\left(\frac{z-m}{de}\right)\left(\frac{z-n}{df}\right) \tag{8}$$

The left side of (8) is a square number, and any pair of the three items on the right side is relatively prime: $\gcd\left(\frac{z-m}{de}, \frac{z-n}{df}\right) = 1$, $\gcd\left(\frac{z}{ef}, \frac{z-m}{de}\right) = 1$, $\gcd\left(\frac{z}{ef} \cdot, \frac{z-n}{df}\right) = 1$. Therefore, each item on the right side of (8) must be a square number:

$$\begin{cases} \frac{z}{ef} = u^2 \\ \frac{z-m}{de} = v^2 \\ \frac{z-n}{df} = w^2 \end{cases} \quad u, v, w \in Z \tag{9}$$

From the first equation in (9), we get $z = efu^2$, and substitute it into the other two equations in (9):

$$\begin{cases} efu^2 - m = dev^2 \\ efu^2 - n = dfw^2 \end{cases}$$

Because $e|m$, $f|n$, we get:

$$\begin{cases} fu^2 - dv^2 = \frac{m}{e} \\ eu^2 - dw^2 = \frac{n}{f} \end{cases} \tag{10}$$

(10) is simutaneous Pell equations. To simplify it, convert it to:

$$\begin{cases} (fu)^2 - fdv^2 = \frac{fm}{e} \\ (eu)^2 - edw^2 = \frac{en}{f} \end{cases} \tag{11}$$

Let $g=fu$, $h=eu$, we have:

$$\begin{cases} g^2 - fd \cdot v^2 = \frac{f(b-a)}{e} \\ h^2 - ed \cdot w^2 = \frac{e(c-a)}{f} \\ \frac{g}{f} = \frac{h}{e} \ \in Z \end{cases} \tag{12}$$

The first two equations in (12) are generalized Pell equations. The currently existing methods to solve the simutaneous Pell equations (11) or (12) can be applied here, [15], [16], [17], [18], [19], [20]. Only if (11) or (12) has integer solutions can we get the solutions of (4). If we get $(u,v,w)$ as one solution to (11), *z*, *x*, and *y* can be solved as: $z = efu^2$, $x = z + a$, and

$$y = \pm\sqrt{(x-a)(x-b)(x-c)}$$

Therefore, by enumerating the factors of *c*-*b*, *b*-*a* and *c*-*a* (named as *d*, *e*, *f* respectively), the equation (4) can be solved by the procedure presented above, which is convenient for computer programming implementation.

Because the sovability of generalized Pell equation can be definitely determined, [21], [22], [23], [24], [25], a sufficient condition of (4)'s solvability is presented based on the above analysis:

*For $y^2 = (x-a)(x-b)(x-c)$, d, e, f are the factors of c-b, b-a and c-a respectively. If for any combination of d, e and f, the generalized equation $g^2 - fd \cdot v^2 = \frac{f(b-a)}{e}$ or $h^2 - ed \cdot w^2 = \frac{e(c-a)}{f}$ has no integer solution, then $y^2 = (x-a)(x-b)(x-c)$ has no integer solution either.*

## 4 The Solution of $y^2 = (x-a)(x^2 + ax + b)$

In this equation, $x^2 + ax + b$ is irreducible in polynomials with integer coefficients. Let $z=x-a$, we have:

$$y^2 = z \cdot \left(z^2 + 3az + (2a^2 + b)\right) \tag{13}$$

Suppose $d = \gcd(z, z^2 + 3az + (2a^2 + b))$ . Then we have $d|2a^2 + b$. Therefore, all the possible values of $d$ belong to the factor set of $2a^2+b$, which can be enumerated finitely. Transform (13) to seperate the coprime parts of the right side:

$$y^2 = d^2 \cdot \frac{z}{d} \cdot \frac{z^2+3az+(2a^2+b)}{d} \quad (14)$$

and

$$\left(\frac{y}{d}\right)^2 = \frac{z}{d} \cdot \frac{z^2+3az+(2a^2+b)}{d} \quad (15)$$

where $\gcd\left(\frac{z}{d}, \frac{z^2+3az+(2a^2+b)}{d}\right) = 1$ , and the left side of (15) is a square number. Therefore, the two items on the right side of (15) must be both square numbers:

$$\begin{cases} \frac{z}{d} = u^2 \\ \frac{z^2+3az+(2a^2+b)}{d} = v^2 \end{cases} \quad u, v \in Z \quad (16)$$

According to (16), we can traverse all values of $u$ in a pre-estimated range (i.e. by the Baker method or the later improvements), calculate $z$ and check whether the second equation in (16) holds true. This is easy to implement by computer programming.

In the following, (16) is further transformed to get a sufficient condition of (13)'s insolvability. Substitute $z = du^2$ into the second equation in (16):

$$dv^2 = d^2u^4 + 3adu^2 + (2a^2 + b) \quad (17)$$

$$4dv^2 = 4d^2u^4 + 12adu^2 + 8a^2 + 4b = (2du^2 + 3a)^2 + (4b - a^2) \quad (18)$$

Let $w = 2du^2 + 3a$, (18) becomes:

$$w^2 - 4dv^2 = a^2 - 4b \quad (19)$$

which has the form of generalized Pell equation. If (19) has no integer solution, the original equation $y^2 = (x - a)(x^2 + ax + b)$ also has no integer solution. Therefore, we have:

*For* $y^2 = (x - a)(x^2 + ax + b)$ , *if for any* ($2a^2+b$)'*s factor d*, *equation* $w^2 - 4dv^2 = a^2 - 4b$ *has no integer solution, then* $y^2 = (x - a)(x^2 + ax + b)$ *has no integer solution either*.

For the generalized Pell equation (19), there are existing methods to determine the existence of integer solution, which can be applied here, [21], [22], [23], [24], [25].

If for some $d$, (19) has integer solution, it has infinite number of integer solutions ($v$,$w$). The solution of generalized Pell equation can be found by the existing methods, [21], [22], [23], [24], [25]. Because (19) may have infinite number of integer solutions, pre-estimation of solution range is needed, which can be obtained by Baker's method or the later improvements. Then in a corresponding finite range of ($v$,$w$), we check whether $\frac{w-3a}{2d}$ is a square number ( $w = 2du^2 + 3a$ , $u^2 = \frac{w-3a}{2d}$). If so, an integer solution ($x$,$y$) of the original equation is found:

$$x = du^2 + a, y = \pm\sqrt{(x - a)(x^2 + ax + b)} \quad (20)$$

Check all the factor $d$ of $2a^2+b$, the integer solution of the original equation can be totally searched, which is convenient for computer programming implementation.

# 5 Conclusion

In this paper, the factorizable elliptic curve equation is converted to limited numbers of simultaneous Pell equations or generalized Pell equations, which can be solved by existing methods. Then the integer solution of the original elliptic curve equation can be effectively searched. The proposed methods can improve the search efficiency in a large range of possible solution. At the same time, two sufficient conditions of insolvability of the factorizable elliptic curve equation are presented. Future study may improve these sufficient conditions with the details in solving generalized Pell equations. Detailed congruence analysis may be introduced into the methods to further improve efficiency.

**Contribution of Individual Authors to the Creation of a Scientific Article (Ghostwriting Policy)**

Xiaodong Zhuang designed the methods and prepared the draft of the manuscript.

Nikos E. Mastorakis did the analysis of the methods and revised the manuscript.

**Sources of Funding for Research Presented in a Scientific Article or Scientific Article Itself**

No funding was received for conducting this study.

**Conflict of Interest**

There is no conflict of interest.